# SMALL TIME PATH BEHAVIOR OF DOUBLE STOCHASTIC INTEGRALS AND APPLICATIONS TO STOCHASTIC CONTROL

By Patrick Cheridito,[1] H. Mete Soner[2] and Nizar Touzi

*Princeton University, Koç University and Centre de Recherche en Economie et Statistique*

We study the small time path behavior of double stochastic integrals of the form $\int_0^t (\int_0^r b(u)\, dW(u))^T\, dW(r)$, where $W$ is a $d$-dimensional Brownian motion and $b$ is an integrable progressively measurable stochastic process taking values in the set of $d \times d$-matrices. We prove a law of the iterated logarithm that holds for all bounded progressively measurable $b$ and give additional results under continuity assumptions on $b$. As an application, we discuss a stochastic control problem that arises in the study of the super-replication of a contingent claim under gamma constraints.

**1. Introduction.** In this paper we study the small time path behavior of double stochastic integrals of the form $V^b(t) = \int_0^t (\int_0^r b(u)\, dW(u))^T\, dW(r)$, where $W$ is a $d$-dimensional Brownian motion and $b$ is an integrable progressively measurable stochastic process taking values in the set of $d \times d$-matrices. We first prove a law of the iterated logarithm under general assumptions. Then, we prove additional results under continuity assumptions on $b$. The results for $V^b$ can easily be generalized to double stochastic integrals of the form $\int_0^t (\int_0^r b(u)\, dM(u))^T\, dM(r)$, for $d$-dimensional martingales $M(t) = \int_0^t m(r)\, dW(r)$ corresponding to regular enough matrix-valued processes $m$.

Results on the small time path behavior of stochastic integrals can be applied to characterize the set of all starting points from which a given controlled continuous-time stochastic process can be driven into a target set

Received March 2004; revised February 2005.
[1]Supported by the Centre de Recherche en Economie et Statistique and the Swiss National Science Foundation.
[2]Supported in part by the Turkish Academy of Sciences.

at a prespecified future time. It is shown in [7] and [8] that under suitable conditions, the set of initial data from which a controlled state process can be steered into a target set, satisfies a dynamic programming principle (DPP), from which a dynamic programming equation (DPE) can be derived. Since in [7] and [8] the control process is constrained to take values in a subset of $\mathbb{R}^d$, the essential step in the derivation of the DPE from the DPP is an analysis of the small time behavior of single stochastic integrals. In [6], the problem of super-replicating a contingent claim under gamma constraints is studied. This problem naturally leads to an analysis of the small time behavior of double stochastic integrals. The results in [6] are obtained under the assumptions that the contingent claim depends on only one underlying asset and that the gamma constraint is an upper bound. In this paper we provide a more extensive study of the small time path behavior of double stochastic integrals than in [6] and discuss the super-replication problem under upper and lower gamma constraints.

In Section 2 we establish the notation and discuss basic examples of double stochastic integrals. The main results of the paper are stated and proved in Section 3, and in Section 4 it is shown how they can be used to find the super-replication price of a contingent claim in the presence of gamma constraints that are more general than in [6]. We keep the presentation in Section 4 simple by making strong assumptions. For a general treatment of the super-replication problem under gamma constraints in a multidimensional framework, we refer the reader to the accompanying paper [2].

**2. Problem formulation and notation.** Let $(\Omega, \mathcal{F}, P)$ be a complete probability space endowed with a filtration $\mathbb{F} := \{\mathcal{F}(t), t \geq 0\}$ that satisfies the usual conditions. We are interested in the small time behavior of double stochastic integrals of the form

$$(2.1) \qquad V^b(t) := \int_0^t \left( \int_0^r b(u) \, dW(u) \right)^T dW(r), \qquad t \geq 0,$$

where $\{W(t), t \geq 0\}$ is a $d$-dimensional Brownian motion on the filtered probability space $(\Omega, \mathcal{F}, \mathbb{F}, P)$, $\{b(t), t \geq 0\}$ is an integrable $\mathbb{F}$-progressively measurable stochastic process with values in $\mathcal{M}^d$, the set of $d \times d$-matrices with real components, and $^T$ denotes the transposition of matrices.

In the easy case where $\{W(t), t \geq 0\}$ is a one-dimensional Brownian motion and $b(t) = \beta$, $t \geq 0$, for some $\beta \in \mathbb{R}$, we have

$$V^b(t) = \frac{\beta}{2}(W^2(t) - t), \qquad t \geq 0.$$

It follows from the law of the iterated logarithm for Brownian motion that

$$(2.2) \qquad \limsup_{t \searrow 0} \frac{2V^\beta(t)}{h(t)} = \beta \qquad \text{for every } \beta \geq 0,$$

where

$$h(t) := 2t \log\log \frac{1}{t}, \qquad t > 0,$$

and the equality in (2.2) is, as all other equalities and inequalities between random variables in this paper, understood in the almost sure sense. On the other hand, it can be deduced from the fact that almost all paths of a one-dimensional Brownian motion cross zero on all time intervals $(0, \varepsilon]$, $\varepsilon > 0$, that

(2.3) $$\limsup_{t \searrow 0} \frac{2V^\beta(t)}{t} = -\beta \qquad \text{for every } \beta < 0.$$

The purpose of this paper is to derive formulae similar to (2.2) and (2.3) when $W = \{W(t), t \geq 0\}$ is a $d$-dimensional Brownian motion and $b = \{b(t), t \geq 0\}$ is a progressively measurable matrix-valued stochastic process. Note that if $b(t) = \beta$, $t \geq 0$, for some fixed symmetric matrix $\beta$, then

$$2V^b(t) = W(t)^T \beta W(t) - \text{Tr}[\beta]t, \qquad t \geq 0,$$

where Tr denotes the trace of a matrix. It is already not completely obvious if the formulae (2.2) and (2.3) have analogs in this case and how they look. In Section 3 we will prove extensions of (2.2) and (2.3) for processes of the form (2.1).

By $I_d$ we denote the $d \times d$ identity matrix. For $y \in \mathbb{R}^n$, we set $|y| := (y_1^2 + \cdots + y_n^2)^{1/2}$, and for $\beta \in \mathcal{M}^d$,

$$|\beta| := \sup_{y \in \mathbb{R}^d, |y|=1} |\beta y|.$$

By $\mathcal{S}^d$ we denote the collection of all symmetric matrices of $\mathcal{M}^d$, and for all $\beta \in \mathcal{S}^d$, we set

$$\lambda_*(\beta) := \min\{y^T \beta y : y \in \mathbb{R}^d, |y| = 1\},$$
$$\lambda^*(\beta) := \max\{y^T \beta y : y \in \mathbb{R}^d, |y| = 1\}.$$

Note that $\lambda^*$ and $\lambda_*$ are continuous, and therefore, measurable functions from $\mathcal{S}^d$ to $\mathbb{R}$. We endow the set $\mathcal{S}^d$ with the usual partial order

$$\beta \geq \alpha \quad \text{if and only if} \quad \lambda_*(\beta - \alpha) \geq 0,$$

and we set $\mathcal{S}^d_+ := \{\beta \in \mathcal{S}^d : \beta \geq 0\}$.

**3. Small time path behavior of double stochastic integrals.** The main results of this section are Theorems 3.1 and 3.3. Corollaries 3.7 and 3.8 are consequences of Theorems 3.1 and 3.3, respectively. Proposition 3.9, whose proof is straightforward, is given because, along with Corollaries 3.7 and 3.8, it is needed in Section 4 of this paper and in the accompanying paper [2].

THEOREM 3.1. *(a) Let $\{b(t), t \geq 0\}$ be an $\mathcal{M}^d$-valued, $\mathbb{F}$-progressively measurable process such that $|b(t)| \leq 1$ for all $t \geq 0$. Then*

$$\limsup_{t \searrow 0} \frac{|2V^b(t)|}{h(t)} \leq 1.$$

*(b) Let $\beta \in \mathcal{S}^d$ with largest eigenvalue $\lambda^*(\beta) \geq 0$. If $\{b(t), t \geq 0\}$ is a bounded, $\mathcal{S}^d$-valued, $\mathbb{F}$-progressively measurable process such that $b(t) \geq \beta$ for all $t \geq 0$, then*

$$\limsup_{t \searrow 0} \frac{2V^b(t)}{h(t)} \geq \lambda^*(\beta).$$

For the proof of Theorem 3.1(a) we need the following exponential estimate:

LEMMA 3.2. *Let $\lambda$ and $T$ be two positive parameters with $2\lambda T < 1$ and $\{b(t), t \geq 0\}$ an $\mathcal{M}^d$-valued, $\mathbb{F}$-progressively measurable process such that $|b(t)| \leq 1$, for all $t \geq 0$. Then*

$$\mathrm{E}[\exp(2\lambda V^b(T))] \leq \mathrm{E}[\exp(2\lambda V^{I_d}(T))].$$

PROOF. We prove this lemma with a standard argument from the theory of stochastic control. We define the processes

$$Y^b(r) := Y(0) + \int_0^r b(u) \, dW(u)$$

and

$$Z^b(t) := Z(0) + \int_0^t (Y^b(r))^T \, dW(r), \qquad t \geq 0,$$

where $Y(0) \in \mathbb{R}^d$ and $Z(0) \in \mathbb{R}$ are some given initial data. Observe that $V^b(t) = Z^b(t)$ when $Y(0) = 0$ and $Z(0) = 0$. We split the argument into three steps.

*Step* 1. It can easily be checked that

(3.1) $$\mathrm{E}[\exp(2\lambda Z^{I_d}(T))|\mathcal{F}(t)] = f(t, Y^{I_d}(t), Z^{I_d}(t)),$$

where, for $t \in [0, T]$, $y \in \mathbb{R}^d$ and $z \in \mathbb{R}$, the function $f$ is given by

$$f(t, y, z) := \mathrm{E}\bigg[\exp\bigg(2\lambda\bigg\{z + \int_t^T (y + W(r) - W(t))^T \, dW(r)\bigg\}\bigg)\bigg]$$

$$= \exp(2\lambda z) \, \mathrm{E}[\exp(\lambda\{2y^T W(T-t) + |W(T-t)|^2 - d(T-t)\})]$$

$$= \mu^{d/2} \exp[2\lambda z - d\lambda(T-t) + 2\mu\lambda^2(T-t)|y|^2],$$

and $\mu := [1 - 2\lambda(T-t)]^{-1}$. Observe that the function $f$ is strictly convex in $y$ and

$$(3.2) \qquad D^2_{yz}f(t,y,z) := \frac{\partial^2 f}{\partial y \, \partial z}(t,y,z) = g^2(t,y,z)y,$$

where $g^2 := 8\mu\lambda^3(T-t)f$ is a positive function of $(t,y,z)$.

*Step* 2. For a matrix $\beta \in \mathcal{M}^d$, we denote by $\mathcal{L}^\beta$ the Dynkin operator associated to the process $(Y^\beta, Z^\beta)$, that is,

$$\mathcal{L}^\beta := D_t + \tfrac{1}{2}\operatorname{Tr}[\beta\beta^T D^2_{yy}] + \tfrac{1}{2}|y|^2 D^2_{zz} + (\beta y)^T D^2_{yz},$$

where $D_\cdot$ and $D^2_{\cdot\cdot}$ denote the gradient and the Hessian operators with respect to the indexed variables. In this step, we intend to prove that for all $t \in [0,T]$, $y \in \mathbb{R}^d$ and $z \in \mathbb{R}$,

$$(3.3) \qquad \max_{\beta \in \mathcal{M}^d, |\beta| \leq 1} \mathcal{L}^\beta f(t,y,z) = \mathcal{L}^{I_d} f(t,y,z) = 0.$$

The second equality can be derived from the fact that the process

$$f(t, Y^{I_d}(t), Z^{I_d}(t)), \qquad t \in [0,T],$$

is a martingale, which can easily be seen from (3.1). The first equality follows from the following two observations: First, note that for each $\beta \in \mathcal{M}^d$ such that $|\beta| \leq 1$, the matrix $I_d - \beta\beta^T$ is in $\mathcal{S}^d_+$. Therefore, there exists a $\gamma \in \mathcal{S}^d_+$ such that

$$I_d - \beta\beta^T = \gamma^2.$$

Since $f$ is convex in $y$, the Hessian matrix $D^2_{yy}f$ is also in $\mathcal{S}^d_+$. It follows that $\gamma D^2_{yy}f(t,x,y)\gamma \in \mathcal{S}^d_+$, and therefore,

$$(3.4) \qquad \begin{aligned} \operatorname{Tr}[D^2_{yy}f(t,x,y)] &- \operatorname{Tr}[\beta\beta^T D^2_{yy}f(t,x,y)] \\ &= \operatorname{Tr}[(I_d - \beta\beta^T)D^2_{yy}f(t,x,y)] \\ &= \operatorname{Tr}[\gamma D^2_{yy}f(t,x,y)\gamma] \geq 0. \end{aligned}$$

Second, it follows from (3.2) and the Cauchy–Schwarz inequality that, for all $\beta \in \mathcal{M}^d$ such that $|\beta| \leq 1$,

$$(3.5) \qquad \begin{aligned} (\beta y)^T D^2_{yz}f(t,y,z) &= g^2(t,y,z)(\beta y)^T y \\ &\leq g^2(t,y,z)|y|^2 \\ &= y^T D^2_{yz}f(t,y,z). \end{aligned}$$

Together, (3.4) and (3.5) imply the first equality in (3.3).

*Step* 3. Let $\{b(t), t \geq 0\}$ be an $\mathcal{M}^d$-valued, $\mathbb{F}$-progressively measurable process such that $|b(t)| \leq 1$ for all $t \geq 0$. We define the sequence of stopping times
$$\tau_k := T \wedge \inf\{t \geq 0 : |Y^b(t)| + |Z^b(t)| \geq k\}, \qquad k \in \mathbb{N}.$$
It follows from Itô's lemma and (3.3) that
$$f(0, Y(0), Z(0)) = f(\tau_k, Y^b(\tau_k), Z^b(\tau_k)) - \int_0^{\tau_k} \mathcal{L}^{b(t)} f(t, Y^b(t), Z^b(t))\, dt$$
$$- \int_0^{\tau_k} [(D_y f)^T b + (D_z f) y^T](t, Y^b(t), Z^b(t))\, dW(t)$$
$$\geq f(\tau_k, Y^b(\tau_k), Z^b(\tau_k))$$
$$- \int_0^{\tau_k} [(D_y f)^T b + (D_z f) y^T](t, Y^b(t), Z^b(t))\, dW(t).$$
Taking expected values and sending $k$ to infinity, we get by Fatou's lemma,
$$\mathrm{E}[\exp(2\lambda Z^{I_d}(T))] = f(0, Y(0), Z(0))$$
$$\geq \liminf_{k \to \infty} \mathrm{E}[f(\tau_k, Y^b(\tau_k), Z^b(\tau_k))]$$
$$\geq \mathrm{E}[f(T, Y^b(T), Z^b(T))]$$
$$= \mathrm{E}[\exp(2\lambda Z^b(T))],$$
which proves the lemma. $\square$

PROOF OF THEOREM 3.1. (a) Let $T > 0$ and $\lambda > 0$ be such that $2\lambda T < 1$. It follows from Doob's maximal inequality for submartingales and Lemma 3.2 that for all $\alpha \geq 0$,
$$P\left[\sup_{0 \leq t \leq T} 2V^b(t) \geq \alpha\right] = P\left[\sup_{0 \leq t \leq T} \exp(2\lambda V^b(t)) \geq \exp(\lambda\alpha)\right]$$
$$\leq \exp(-\lambda\alpha)\, \mathrm{E}[\exp(2\lambda V^b(T))]$$
(3.6)
$$\leq \exp(-\lambda\alpha)\, \mathrm{E}[\exp(2\lambda V^{I_d}(T))]$$
$$= \exp(-\lambda\alpha) \exp(-\lambda d T)(1 - 2\lambda T)^{-d/2}.$$
Now, take $\theta, \eta \in (0, 1)$, and set for all $k \in \mathbb{N}$,
$$\alpha_k := (1 + \eta)^2 h(\theta^k) \quad \text{and} \quad \lambda_k := [2\theta^k(1 + \eta)]^{-1}.$$
It follows from (3.6) that for all $k \in \mathbb{N}$,
$$P\left[\sup_{0 \leq t \leq \theta^k} 2V^b(t) \geq (1 + \eta)^2 h(\theta^k)\right]$$
$$\leq \exp\left(-\frac{d}{2(1 + \eta)}\right) \left(1 + \frac{1}{\eta}\right)^{d/2} \left(k \log \frac{1}{\theta}\right)^{-(1+\eta)}.$$

Since
$$\sum_{k=1}^{\infty} k^{-(1+\eta)} < \infty,$$

it follows from the Borel–Cantelli lemma that, for almost all $\omega \in \Omega$, there exists a natural number $K^{\theta,\eta}(\omega)$ such that for all $k \geq K^{\theta,\eta}(\omega)$,
$$\sup_{0 \leq t \leq \theta^k} 2V^b(t,\omega) < (1+\eta)^2 h(\theta^k).$$

In particular, for all $t \in (\theta^{k+1}, \theta^k]$,
$$2V^b(t,\omega) < (1+\eta)^2 h(\theta^k) \leq (1+\eta)^2 \frac{h(t)}{\theta}.$$

Hence,
$$\limsup_{t \searrow 0} \frac{2V^b(t)}{h(t)} \leq \frac{(1+\eta)^2}{\theta}.$$

By letting $\theta$ tend to 1 and $\eta$ to zero along the rationals, we conclude that
$$\limsup_{t \searrow 0} \frac{2V^b(t)}{h(t)} \leq 1.$$

On the other hand,
$$\liminf_{t \searrow 0} \frac{2V^b(t)}{h(t)} = -\limsup_{t \searrow 0} \frac{2V^{-b}(t)}{h(t)} \geq -1,$$

and the proof of part (a) is complete.

(b) There exists a constant $c > 0$ such that for all $t \geq 0$,

(3.7) $$cI_d \geq b(t) \geq \beta \geq -cI_d,$$

and an orthogonal $d \times d$-matrix $U$ such that
$$\tilde{\beta} := U\beta U^T = \operatorname{diag}[\lambda^*(\beta), \lambda_2, \ldots, \lambda_d],$$

where $\lambda^*(\beta) \geq \lambda_2 \geq \cdots \geq \lambda_d$ are the ordered eigenvalues of the matrix $\beta$. Let
$$\tilde{\gamma} := \operatorname{diag}[3c, c, \ldots, c] \quad \text{and} \quad \gamma := U^T \tilde{\gamma} U.$$

It follows from (3.7) that for all $t \geq 0$,
$$\gamma - \beta \geq \gamma - b(t) \geq 0,$$

which implies that
$$|\gamma - b(t)| \leq |\gamma - \beta| = \lambda^*(\gamma - \beta) = \lambda^*(\tilde{\gamma} - \tilde{\beta}) = 3c - \lambda^*(\beta).$$

Hence, by part (a),

$$\limsup_{t \searrow 0} \frac{2V^{\gamma-b}(t)}{h(t)} \leq 3c - \lambda^*(\beta). \tag{3.8}$$

Note that the transformed Brownian motion,

$$\tilde{W}(t) := UW(t), \qquad t \geq 0,$$

is again a $d$-dimensional Brownian motion and

$$\begin{aligned}
\limsup_{t \searrow 0} \frac{2V^{\gamma}(t)}{h(t)} &= \limsup_{t \searrow 0} \frac{W(t)^T \gamma W(t) - \text{Tr}(\gamma)t}{h(t)} \\
&= \limsup_{t \searrow 0} \frac{\tilde{W}(t)^T \tilde{\gamma} \tilde{W}(t) - \text{Tr}(\gamma)t}{h(t)} \\
&= \limsup_{t \searrow 0} \frac{\tilde{W}(t)^T \tilde{\gamma} \tilde{W}(t)}{h(t)} \\
&\geq \limsup_{t \searrow 0} 3c \frac{(\tilde{W}_1(t))^2}{h(t)} = 3c.
\end{aligned} \tag{3.9}$$

It follows from (3.9) and (3.8) that

$$\limsup_{t \searrow 0} \frac{2V^b(t)}{h(t)} \geq \limsup_{t \searrow 0} \frac{2V^{\gamma}(t)}{h(t)} - \limsup_{t \searrow 0} \frac{2V^{\gamma-b}(t)}{h(t)}$$
$$\geq 3c - (3c - \lambda^*(\beta)) = \lambda^*(\beta),$$

which proves part (b) of the theorem. □

In the next theorem we refine the statements of Theorem 3.1 under stronger assumptions.

THEOREM 3.3. *Let $\{b(t), t \geq 0\}$ be an $\mathcal{M}^d$-valued, $\mathbb{F}$-progressively measurable process such that*

$$\int_0^t |b(r)|^2 \, dr < \infty \qquad \text{for all } t \geq 0.$$

*Assume that $b(0)$ is a deterministic element of $\mathcal{S}^d$, and there exists a random variable $\varepsilon > 0$ such that almost surely,*

$$\int_0^t |b(r) - b(0)|^2 \, dr = O(t^{1+\varepsilon}) \qquad \text{for } t \searrow 0. \tag{3.10}$$

(a) *If $\lambda^*(b(0)) \leq 0$, then*

$$\limsup_{t \searrow 0} \frac{2V^b(t)}{t} = -\text{Tr}[b(0)].$$

(b) *If* $\lambda^*(b(0)) \geq 0$, *then*

$$\limsup_{t \searrow 0} \frac{2V^b(t)}{h(t)} = \lambda^*(b(0)).$$

REMARK 3.4. Note that for deterministic $\varepsilon > 0$, condition (3.10) follows if there exists a constant $C > 0$ such that

(3.11) $$\mathrm{E}[|b(t) - b(0)|^2] \leq Ct^{2\varepsilon} \qquad \text{for } t \geq 0.$$

Indeed, if (3.11) is satisfied, then

$$\mathrm{E}\left[\int_0^1 \frac{|b(r) - b(0)|^2}{r^{1+\varepsilon}} \, dr\right] < \infty,$$

therefore,

$$\int_0^1 \frac{|b(r) - b(0)|^2}{r^{1+\varepsilon}} \, dr < \infty,$$

and we have for all $t \in [0,1]$,

$$\int_0^t |b(r) - b(0)|^2 \, dr \leq \int_0^t \frac{|b(r) - b(0)|^2}{r^{1+\varepsilon}} \, dr \, t^{1+\varepsilon} \leq \int_0^1 \frac{|b(r) - b(0)|^2}{r^{1+\varepsilon}} \, dr \, t^{1+\varepsilon}.$$

To prove Theorem 3.3 we need the following.

LEMMA 3.5. *Let* $\{W(t), t \geq 0\}$ *be a d-dimensional Brownian motion and* $\beta \in \mathcal{M}^d$. *Then*

(3.12) $$\liminf_{t \searrow 0} \frac{1}{t}|W(t)^T \beta W(t)| = 0.$$

PROOF. It follows from the self-similarity property of $\{W(t), t \geq 0\}$ that the Gaussian sequence,

$$X(n) := e^{n/2} W(e^{-n}), \qquad n \in \mathbb{Z},$$

is stationary, and the fact that

$$\lim_{n \to \infty} \mathrm{E}[X(n)^T X(0)] = 0$$

implies that it is ergodic (see, e.g., Section V.3 in [5]). Hence, the sequence

$$Y(n) := |X(n)^T \beta X(n)| = e^n |W(e^{-n})^T \beta W(e^{-n})|, \qquad n \in \mathbb{Z},$$

is stationary and ergodic as well. Therefore, we can apply the ergodic theorem (see, e.g., Theorem V.3.3 in [5]) to conclude that for all $\delta > 0$,

$$\lim_{n \to \infty} \frac{1}{n} \sum_{j=0}^{n-1} \mathbf{1}_{[0,\delta]}(Y(j)) = \mathrm{E}[\mathbf{1}_{[0,\delta]}(Y(0))] = P[Y(0) \leq \delta] > 0.$$

This shows that
$$\liminf_{n \to \infty} Y(n) = 0,$$
which implies (3.12). □

PROOF OF THEOREM 3.3. Since $b(0)$ is symmetric, we have for all $t \geq 0$,

(3.13)  $2V^b(t) = 2V^{b(0)}(t) + 2V^{\tilde{b}}(t) = W(t)^T b(0) W(t) - \text{Tr}[b(0)]t + 2V^{\tilde{b}}(t),$

where
$$\tilde{b}(t) := b(t) - b(0), \qquad t \geq 0.$$

Denote by $M_j$ the $j$th component of the $d$-dimensional local martingale $\int_0^r \tilde{b}(u) \, dW(u)$, $r \geq 0$. It follows from assumption (3.10) that the quadratic variation process $\langle M_j \rangle$ satisfies almost surely,

$$\langle M_j \rangle(r) = \int_0^r \sum_{k=1}^d \tilde{b}_{jk}^2(u) \, du = O(r^{1+\varepsilon}) \qquad \text{for } t \searrow 0.$$

By the Dambis–Dubins–Schwarz theorem (see, e.g., Section V.1 in [4]), there exists a Brownian motion $B_j$ such that $M_j(r) = B_j \circ \langle M_j \rangle(r)$, $r \geq 0$. Hence, it follows from the law of the iterated logarithm for Brownian motion that almost surely,

$$M_j^2(r) = O(r^{1+\varepsilon/2}) \qquad \text{for } r \searrow 0.$$

This implies that almost surely,

$$\langle V^{\tilde{b}} \rangle(t) = \int_0^t \sum_{j=1}^d M_j^2(r) \, dr = O(t^{2+\varepsilon/2}) \qquad \text{for } t \searrow 0,$$

and another application of the Dambis–Dubins–Schwarz theorem yields

(3.14) $$\lim_{t \searrow 0} \frac{V^{\tilde{b}}(t)}{t} = 0.$$

(a) If $\lambda^*(b(0)) \leq 0$, then for all $t \geq 0$,
$$W(t)^T b(0) W(t) \leq 0,$$
and part (a) of the theorem can be deduced from (3.13), (3.14) and Lemma 3.5.

(b) If $\lambda^*(b(0)) \geq 0$, it follows from Theorem 3.1(b) that

(3.15) $$\limsup_{t \searrow 0} \frac{2V^{b(0)}(t)}{h(t)} \geq \lambda^*(b(0)).$$

To show that actually,

(3.16) $$\limsup_{t \searrow 0} \frac{2V^{b(0)}(t)}{h(t)} = \lambda^*(b(0)),$$

we denote by $\lambda^*(b(0)) = \lambda_1 \geq \lambda_2 \geq \cdots \geq \lambda_d$ the ordered eigenvalues of $b(0)$. There exists a positive integer $k \leq d$ such that $\lambda_1 \geq \cdots \geq \lambda_k \geq 0$ and, in case that $k < d$, $0 > \lambda_{k+1} \geq \cdots \geq \lambda_d$. Let $U$ be an orthogonal $d \times d$-matrix such that

$$Ub(0)U^T = \text{diag}[\lambda_1, \lambda_2, \ldots, \lambda_d].$$

The process

$$\tilde{W}(t) := UW(t), \qquad t \geq 0,$$

is again a $d$-dimensional Brownian motion, and for all $t \geq 0$,

$$\limsup_{t \searrow 0} \frac{2V^{b(0)}(t)}{h(t)} = \limsup_{t \searrow 0} \frac{W(t)^T b(0) W(t) - \text{Tr}[b(0)]t}{h(t)}$$
$$= \limsup_{t \searrow 0} \frac{\sum_{j=1}^d \lambda_j (\tilde{W}_j(t))^2}{h(t)} \leq \limsup_{t \searrow 0} \frac{\sum_{j=1}^k \lambda_j (\tilde{W}_j(t))^2}{h(t)}$$
$$\leq \lambda_1 = \lambda^*(b(0)),$$

where the last inequality follows from Theorem 3.1(a). This and (3.15) imply (3.16), which, along with (3.13) and (3.14), proves part (b) of the theorem. □

Our proof of Theorem 3.3 is based on the decomposition (3.13) and the estimate (3.14). The next example shows that (3.14) need no longer be true if assumption (3.10) is replaced by the condition that almost surely,

$$|b(t) - b(0)| \to 0 \qquad \text{as } t \to 0.$$

Whether Theorem 3.3, or some variant of it, can be proved under weaker assumptions is an open question.

EXAMPLE 3.6. Let $d = 1$ and $b(t) = 1/\log\log\log(1/t)$. Then,

$$\int_0^t \int_0^r b(u)\, dW(u)\, dW(r) = W(t) \int_0^t b(r)\, dW(r)$$
$$- \int_0^t b(r) W(r)\, dW(r) - \int_0^t b(r)\, dr$$
$$= W(t)\left[W(t)b(t) - \int_0^t W(r)\, db(r)\right] - \int_0^t b(r)\, dr$$

(3.17)
$$-\tfrac{1}{2}\left[W^2(t)b(t) - \int_0^t W^2(r)\,db(r) - \int_0^t b(r)\,dr\right]$$
$$= \tfrac{1}{2}W^2(t)b(t) - W(t)\int_0^t W(r)\,db(r)$$
$$+ \tfrac{1}{2}\int_0^t W^2(r)\,db(r) - \tfrac{1}{2}\int_0^t b(r)\,dr.$$

Clearly, $\int_0^t b(r)\,dr = o(t)$, as $t \searrow 0$. Since
$$b'(r) = \frac{1}{r}\frac{1}{\log(1/r)}\frac{1}{\log\log(1/r)}\left(\frac{1}{\log\log\log(1/r)}\right)^2,$$
it follows from the law of the iterated logarithm for Brownian motion that for $t \searrow 0$,
$$\int_0^t W^2(r)\,db(r) = \int_0^t W^2(r)b'(r)\,dr$$
$$= O\left(\int_0^t r\log\log\frac{1}{r}b'(r)\,dr\right)$$
$$= o\left(\int_0^t 1\,dr\right) = o(t).$$

Similarly, for $t \searrow 0$,
$$W(t)\int_0^t W(r)\,db(r) = W(t)\int_0^t W(r)b'(r)\,dr$$
$$= O\left(\sqrt{t\log\log\frac{1}{t}}\int_0^t \sqrt{r\log\log\frac{1}{r}}b'(r)\,dr\right)$$
$$= o\left(\sqrt{t\log\log\frac{1}{t}}\int_0^t r^{-1/2}\,dr\right)$$
$$= o\left(t\sqrt{\log\log\frac{1}{t}}\right).$$

Since
$$t\sqrt{\log\log\frac{1}{t}} = o\left(t\frac{\log\log(1/t)}{\log\log\log(1/t)}\right) \qquad \text{as } t \searrow 0,$$
it follows from (3.17) that
$$\limsup_{t \searrow 0} \frac{\int_0^t \int_0^r b(u)\,dW(u)\,dW(r)}{t\log\log(1/t)/(\log\log\log(1/t))}$$

$$= \limsup_{t \searrow 0} \frac{(1/2)W^2(t)b(t)}{t \log\log(1/t)/(\log\log\log(1/t))} = 1.$$

The next two corollaries are straightforward consequences of Theorems 3.1 and 3.3, respectively.

COROLLARY 3.7. *Let $\{M(t), t \geq 0\}$ be an $\mathbb{R}^d$-valued martingale defined by*

$$M(t) := \int_0^t m(r)\,dW(r), \qquad t \geq 0,$$

*where $\{m(t), t \geq 0\}$ is an $\mathcal{M}^d$-valued, $\mathbb{F}$-progressively measurable process such that*

$$\int_0^t |m(r)|^2\,dr < \infty \qquad \text{for all } t \geq 0,$$

*and there exists a random variable $\varepsilon > 0$ so that almost surely,*

(3.18) $$\int_0^t |m(r) - m(0)|^2\,dr = O(t^{1+\varepsilon}) \qquad \text{for } t \searrow 0.$$

(a) *Let $\{b(t), t \geq 0\}$ be a bounded $\mathcal{M}^d$-valued, $\mathbb{F}$-progressively measurable process such that for all $t \geq 0$, $|m(0)^T b(t) m(0)| \leq 1$. Then*

$$\limsup_{t \searrow 0} \frac{2}{h(t)} \left| \int_0^t \left( \int_0^u b(u)\,dM(u) \right)^T dM(r) \right| \leq 1.$$

(b) *Let $\beta$ be a bounded, $\mathcal{F}(0)$-measurable, $\mathcal{S}^d$-valued random variable with $\lambda^*(\beta) \geq 0$. If $\{b(t), t \geq 0\}$ is a bounded, $\mathcal{S}^d$-valued, $\mathbb{F}$-progressively measurable process such that for all $t \geq 0$,*

$$m(0)^T b(t) m(0) \geq \beta,$$

*then*

$$\limsup_{t \searrow 0} \frac{2}{h(t)} \int_0^t \left( \int_0^r b(u)\,dM(u) \right)^T M(r) \geq \lambda^*(\beta).$$

PROOF. It can easily be checked that

$$\int_0^t \left( \int_0^r b(u) m(u)\,dW(u) \right)^T m(r)\,dW(r)$$
$$= \int_0^t \left( \int_0^r c(u)\,dW(u) \right)^T W(r) + R_1(t) + R_2(t),$$

where
$$c(t) := m(0)^T b(t) m(0),$$
$$R_1(t) := \int_0^t \left( \int_0^r b(u)[m(u) - m(0)] \, dW(u) \right)^T m(0) \, dW(r),$$
$$R_2(t) := \int_0^t \left( \int_0^r b(u) m(u) \, dW(u) \right)^T [m(r) - m(0)] \, dW(r).$$

As in the proof of Theorem 3.3 it can be deduced from assumption (3.18) and the Dambis–Dubins–Schwarz theorem that
$$\lim_{t \searrow 0} \frac{R_1(t)}{t} = \lim_{t \searrow 0} \frac{R_2(t)}{t} = 0.$$

In particular,
$$\lim_{t \searrow 0} \frac{R_1(t)}{h(t)} = \lim_{t \searrow 0} \frac{R_2(t)}{h(t)} = 0.$$

Now, part (a) of the corollary follows from Theorem 3.1(a). Furthermore, by conditioning on $\sigma(\beta)$, we can assume that $\beta$ is deterministic and deduce part (b) of the corollary from Theorem 3.1(b). $\square$

COROLLARY 3.8. *Let $\{M(t), t \geq 0\}$ be an $\mathbb{R}^d$-valued martingale defined by*
$$M(t) = \int_0^t m(r) \, dW(r), \qquad t \geq 0,$$
*where $\{m(t), t \geq 0\}$ is an $\mathcal{M}^d$-valued, $\mathbb{F}$-progressively measurable process such that*
$$\int_0^t |m(r)|^2 \, dr < \infty \qquad \text{for all } t \geq 0.$$
*Let $\{b(t), t \geq 0\}$ be a bounded, $\mathcal{M}^d$-valued, $\mathbb{F}$-progressively measurable process such that $b(0)$ is $\mathcal{S}^d$-valued, and assume there exists a random variable $\varepsilon > 0$ such that almost surely,*
$$\int_0^t |m(r) - m(0)|^2 \, dr = O(t^{1+\varepsilon})$$
*and*
$$\int_0^t |b(r) - b(0)|^2 \, dr = O(t^{1+\varepsilon}) \qquad \text{for } t \searrow 0.$$

(a) *If $\lambda^*(m(0)^T b(0) m(0)) \leq 0$, then*
$$\limsup_{t \searrow 0} \frac{2}{t} \int_0^t \left( \int_0^r b(u) \, dM(u) \right)^T dM(r) = -\operatorname{Tr}[m(0)^T b(0) m(0)].$$

(b) If $\lambda^*(m(0)^T b(0) m(0)) \geq 0$, then

$$\limsup_{t \searrow 0} \frac{2}{h(t)} \int_0^t \left( \int_0^r b(u) \, dM(u) \right)^T dM(r) = \lambda^*(m(0)^T b(0) m(0)).$$

PROOF. As in the proof of Corollary 3.7 we decompose

$$\int_0^t \left( \int_0^r b(u) m(u) \, dW(u) \right)^T m(r) \, dW(r)$$

into

$$\int_0^t \left( \int_0^r c(u) \, dW(u) \right)^T dW(r) + R_1(t) + R_2(t),$$

where

$$c(t) := m(0)^T b(t) m(0),$$

$$\frac{1}{t} R_1(t) := \frac{1}{t} \int_0^t \left( \int_0^r b(u)[m(u) - m(0)] \, dW(u) \right)^T m(0) \, dW(r) \to 0$$

$$\text{for } t \searrow 0,$$

$$\frac{1}{t} R_2(t) := \frac{1}{t} \int_0^t \left( \int_0^r b(u) m(u) \, dW(u) \right)^T [m(r) - m(0)] \, dW(r) \to 0$$

$$\text{for } t \searrow 0.$$

It follows from the assumptions that $c$ satisfies almost surely,

$$\int_0^t |c(r) - c(0)|^2 \, dr = O(t^{1+\varepsilon}) \qquad \text{for } t \searrow 0,$$

and by conditioning on $\sigma(c(0))$, we can assume that $c(0)$ is deterministic. Then, the corollary follows from Theorem 3.3. □

PROPOSITION 3.9. *Let $\{a(t), t \geq 0\}$ and $\{m(t), t \geq 0\}$ be two $\mathbb{F}$-progressively measurable processes taking values in $\mathbb{R}^d$ and $\mathcal{M}^d$, respectively. Assume that $\{a(t), t \geq 0\}$ is bounded,*

$$\int_0^t |m(r)|^2 \, dr < \infty \qquad \text{for all } t \geq 0,$$

*and there exists a $(0,1]$-valued random variable $\varepsilon$ such that almost surely,*

(3.19) $$\int_0^t r^2 |m(r)|^2 \, dr = O(t^{3-\varepsilon}) \qquad \text{for } t \searrow 0.$$

*Then,*

$$\lim_{t \searrow 0} t^{-3/2+\varepsilon} \int_0^t \left( \int_0^r a(u) \, du \right)^T m(r) \, dW(r) = 0.$$

REMARK 3.10. It can easily be shown that

$$\sup_{t\geq 0} \mathrm{E}[|m(t)|^2] < \infty \tag{3.20}$$

implies condition (3.19) for every constant $\varepsilon \in (0,1]$. Indeed, it follows from (3.20) that for every constant $\varepsilon \in (0,1]$,

$$\mathrm{E}\bigg[\int_0^1 \frac{|m(r)|^2}{r^{1-\varepsilon}}\, dr\bigg] < \infty$$

and therefore,

$$\int_0^1 \frac{|m(r)|^2}{r^{1-\varepsilon}}\, dr < \infty.$$

Moreover, for all $t \in [0,1]$,

$$\int_0^t r^2 |m(r)|^2\, dr \leq \int_0^1 \frac{|m(r)|^2}{r^{1-\varepsilon}}\, dr\, t^{3-\varepsilon}.$$

PROOF OF PROPOSITION 3.9. Denote $X(t) = \int_0^t (\int_0^r a(u)\, du)^T m(r)\, dW(r)$, $t \geq 0$. By assumption (3.19), the quadratic variation process $\langle X \rangle$ satisfies almost surely,

$$\langle X \rangle(t) = O(t^{3-\varepsilon}) \qquad \text{for } t \searrow 0.$$

Now, the proposition can be deduced from the Dambis–Dubins–Schwarz theorem. $\square$

**4. Applications to stochastic control.** In this section we show how results on the small time behavior of stochastic integrals can be applied to derive partial differential equations from gamma constraints on hedging strategies. Since these partial differential equations will be derived from a dynamic programming principle (DPP), we refer to them as dynamic programming equations (DPEs).

4.1. *Super-replication under gamma constraints.* For the sake of simplicity of presentation, we here consider a financial market that consists of only two tradable assets. Markets with more assets are considered in the accompanying paper [2]. Let $T > 0$ be a finite time horizon, let $\{W(t), t \in [0,T]\}$ be a one-dimensional Brownian motion and let $\mathbb{F}^W = \{\mathcal{F}^W(t), t \in [0,T]\}$ be the smallest filtration that contains the filtration generated by $\{W(t), t \geq 0\}$ and satisfies the usual conditions. We take the first asset as numéraire and assume that the price of the second one is given by

$$S(r) := S_0 \exp\bigg\{\bigg(\mu - \frac{\sigma^2}{2}\bigg)r + \sigma W(r)\bigg\}, \qquad r \in [0,T],$$

for some constants $S_0 > 0$, $\mu \in \mathbb{R}$ and $\sigma > 0$. By possibly passing to an equivalent probability measure, we can assume that $\mu = 0$. Then, given $S(t) = s$ for some $(t, s) \in [0, T) \times (0, \infty)$, the further evolution of $S$ is

$$(4.1) \qquad S(r) := s \exp\left\{\sigma[W(r) - W(t)] - \frac{\sigma^2}{2}(r - t)\right\}, \qquad r \in [t, T].$$

A self-financing trading strategy that is only based on information coming from observations of the price $\{S(r), r \in [0, T]\}$, can be described by an $\mathbb{F}^W$-progressively measurable process $\{Y(r), r \in [t, T]\}$ that is integrable with respect to $\{S(r), r \in [t, T]\}$ and denotes the number of shares of the second asset held at any given time. Then, the wealth process is given by

$$(4.2) \qquad X(r) = X(t) + \int_t^r Y(u)\, dS(u), \qquad r \in [t, T],$$

and the number of shares of the first asset held at time $r$ is $X(r) - Y(r)S(r)$.

We consider a contingent claim with a time-$T$ payoff given by $g(S(T))$, where $g: (0, \infty) \to [0, \infty)$ is a measurable function such that $g(S(T)) \in L^1(P)$. For the corresponding Black–Scholes hedging strategy $\{Y^{\mathrm{BS}}(r), r \in [t, T]\}$ we have

$$\mathrm{E}[g(S(T))|\mathcal{F}(t)] + \int_t^T Y^{\mathrm{BS}}(r)\, dS(r) = g(S(T)),$$

that is, starting with initial capital $\mathrm{E}[g(S(T))|\mathcal{F}(t)]$ at time $t$, the Black–Scholes strategy replicates the contingent claim. If one requires the hedging strategy to satisfy constraints other than conditions that exclude arbitrage opportunities, one cannot hope that the contingent claim is still replicable, but in many cases, it is possible to super-replicate it with finite initial wealth. A gamma constraint is a restriction on the variation of the hedging strategy due to changes in the underlying asset. To be able to express gamma constraints more explicitly, we require the process $Y$ to be of the form

$$(4.3) \qquad Y(r) = y + \int_t^r \alpha(u)\, du + \int_t^r \gamma(u)\, dS(u), \qquad r \in [t, T],$$

for $y \in \mathbb{R}$ and $\alpha, \gamma$ bounded, $\mathbb{F}^W$-progressively measurable processes. Then, a self-financing trading strategy is determined by the starting value $y$ and a pair of bounded, $\mathbb{F}^W$-progressively measurable processes $\nu = (\alpha, \gamma)$. By a gamma constraint we mean a restriction on the process $\gamma$. In the following we consider gamma constraints of the form:

$$(4.4) \qquad \Gamma_* \leq S^2(r)\gamma(r) \leq \Gamma^*, \qquad r \in [t, T],$$

where $\Gamma_* < \Gamma^*$ are two given constants. We call a control process $\nu = (\alpha, \gamma)$ admissible if $\alpha$ and $\gamma$ are bounded, $\mathbb{F}^W$-progressively measurable processes and $\gamma$ satisfies the constraint (4.4).

To emphasize the dependence on the initial data, we denote by $(S_{t,s}, X^\nu_{t,s,x,y}, Y^\nu_{t,s,y})$ the processes given by (4.1), (4.2) and (4.3) corresponding to the admissible control $\nu$ and the initial data $(S_{t,s}, X^\nu_{t,s,x,y}, Y^\nu_{t,s,y})(t) = (s, x, y)$. The collection of admissible controls $\nu$ is denoted by $\mathcal{A}_{t,s}$. From the boundedness of $\alpha$ and $\gamma$ it can be deduced that for all $\nu \in \mathcal{A}_{t,s}$, $\sup_{t \leq r \leq T} \mathrm{E}[\{Y^\nu_{t,s,y}(r) S_{t,s}(r)\}^2] < \infty$, and therefore, $X^\nu_{t,s,x,y}$ is a square-integrable martingale. In particular, admissible controls do not lead to arbitrage.

The problems

$$(4.5) \quad w(t,s,y) := \inf\{x : X^\nu_{t,s,x,y}(T) \geq g(S_{t,s}(T)) \text{ for some } \nu \in \mathcal{A}\}$$

and

$$\begin{aligned}(4.6) \quad v(t,s) &:= \inf_{y \in \mathbb{R}} w(t,s,y) \\ &= \inf\{x : X^\nu_{t,s,x,y}(T) \geq g(S_{t,s}(T)) \\ &\qquad \text{for some } y \in \mathbb{R} \text{ and some } \nu \in \mathcal{A}\}\end{aligned}$$

can both be viewed as stochastic target problems. Problem (4.5) is very similar to the one treated in [7] and leads to the study of the small time behavior of single stochastic integrals. In problem (4.6), $Y$ is no longer a state variable, and one is naturally led to an analysis of the small time behavior of double stochastic integrals.

In the next two subsections we derive DPEs for $w$ and $v$. Our main objective is to show how one can find these DPEs and where results on the small time behavior of stochastic integrals are needed. To avoid the use of the theory of viscosity solutions and lengthy approximation arguments, we will make some strong assumptions along the way. In particular, we will assume that the infima in (4.5) and (4.6) are attained and the functions $w$ and $v$ are smooth. Also, we will only show that $w$ and $v$ are supersolutions of the corresponding DPEs. A more detailed discussion of the super-replication problem under gamma constraints and rigorous proofs without simplifying assumptions can be found in [2].

4.2. *DPE for the value function w.* We derive the DPE for $w$ in three steps.

*Step* 1: *Dynamic programming principle.* We assume that for each $(t,s,y) \in [0,T) \times (0,\infty) \times \mathbb{R}$, there exists an admissible control $\nu = (\alpha, \gamma)$ such that

$$X^\nu_{t,s,x,y}(T) \geq g(S_{t,s}(T)) \qquad \text{where } x = w(t,s,y).$$

Let $\tau$ be an $\mathbb{F}^W$-stopping time with values in $(t, T]$. For each $\delta > 0$, we define $\tau_\delta := \tau \wedge (t + \delta)$, and we set $(\hat{s}, \hat{x}, \hat{y}) := (S_{t,s}, X^\nu_{t,s,x,y}, Y^\nu_{t,s,y})(\tau_\delta)$. It can be deduced from

$$X^\nu_{\tau_\delta, \hat{s}, \hat{x}, \hat{y}}(T) \geq g(S_{\tau_\delta, \hat{s}}(T))$$

that $\hat{x} \geq w(\tau_\delta, \hat{s}, \hat{y})$, that is,

$$(4.7) \qquad w(t,s,y) + \int_t^{\tau_\delta} Y^\nu_{t,s,y}(r) \, dS_{t,s}(r) \geq w(\tau_\delta, S_{t,s}(\tau_\delta), Y^\nu_{t,s,y}(\tau_\delta)).$$

*Step* 2: *Application of Itô's lemma.* We further assume that the value function $w$ is smooth. Then, we can apply Itô's lemma in (4.7) to get for all $\delta > 0$:

$$(4.8) \qquad -\int_t^{\tau_\delta} \xi(r) \, dr - \int_t^{\tau_\delta} \psi(r) \, dS_{t,s}(r) \geq 0,$$

where

$$\xi(r) := \mathcal{G}^\nu w(r \wedge \tau, S_{t,s}(r \wedge \tau), Y^\nu_{t,s,y}(r \wedge \tau)),$$
$$\psi(r) := (w_s + \gamma w_y)(r \wedge \tau, S_{t,s}(r \wedge \tau), Y^\nu_{t,s,y}(r \wedge \tau)) - Y^\nu_{t,s,y}(r \wedge \tau),$$

and $\mathcal{G}^\nu$ is the Dynkin operator associated to the two-dimensional process $(S, Y^\nu)$:

$$(4.9) \quad \begin{aligned} \mathcal{G}^\nu w(t,s,y) &:= w_t(t,s,y) + \alpha(t) w_y(t,s,y) \\ &\quad + \tfrac{1}{2} \sigma^2 s^2 w_{ss}(t,s,y) + \tfrac{1}{2} \gamma(t)^2 \sigma^2 s^2 w_{yy}(t,s,y) \\ &\quad + \gamma(t) \sigma^2 s^2 w_{sy}(t,s,y). \end{aligned}$$

If we set

$$\tau := \inf\{r \geq t : |Y^\nu_{t,s,y}(r) - y| + |\log S_{t,s}(r) - \log s| > K\} \wedge T,$$

for some constant $K > 0$, then the processes $\xi$ and $\psi$ are bounded.

*Step* 3: *Small time path behavior of single stochastic integrals.* Since $\xi$ is bounded, it follows from (4.8) that there exists a constant $L > 0$ such that

$$(4.10) \quad \int_t^r -\psi(u) \, dS_{t,s}(u) = \int_t^r -\psi(u) S_{t,s}(u) \sigma \, dW(u) \geq -L(r - t)$$

for all $r \in [t, \tau]$.

By the Dambis–Dubins–Schwarz theorem, there exists a Brownian motion $\{B(r), r \geq 0\}$ such that

$$\int_t^r -\psi(u) \, dS_{t,s}(u) = B\left(\int_t^r \psi^2(u) S^2_{t,s}(u) \sigma^2 \, du\right), \qquad r \in [t, T].$$

Hence, it follows from (4.10) and the law of the iterated logarithm for Brownian motion that for all $\varepsilon, \delta > 0$,

$$P[|\psi(u)| \geq \varepsilon \text{ for all } u \in [t, t + \delta]] = 0.$$

By the definition of $\psi$ and the gamma constraint (4.4) on the process $\gamma$, this provides

$$(4.11) \quad -S_\Gamma(-w_y(t,s,y)) \leq s^2(y - w_s(t,s,y)) \leq S_\Gamma(w_y(t,s,y)),$$

where $S_\Gamma$ is the support function of the interval $[\Gamma_*, \Gamma^*]$ defined by

$$S_\Gamma(u) := \sup_{\Gamma_* \leq c \leq \Gamma^*} uc, \quad u \in \mathbb{R}.$$

Since $\psi$ is bounded, we can take expected values in (4.8) and divide by $\delta$ to obtain

$$\mathrm{E}\left[-\frac{1}{\delta}\int_t^{\tau_\delta} \xi(r)\,dr\right] \geq 0,$$

which, in the limit $\delta \to 0$, implies that

$$(4.12) \quad -\mathcal{G}w(t,s,y) := \sup\{-\mathcal{G}^{(a,c)}w(t,s,y) : a \in \mathbb{R} \text{ and } \Gamma_* \leq s^2 c \leq \Gamma^*\} \geq 0,$$

where $\mathcal{G}^{(a,c)}$ is given by (4.9). Combining (4.11) and (4.12), we obtain

$$G(s, y, w_t(t,s,y), Dw(t,s,y), D^2w(t,s,y))$$
$$:= \min\{-\mathcal{G}w(t,s,y); s^2 y - [s^2 w_s - S_\Gamma(-w_y)](t,s,y);$$
$$- s^2 y + [s^2 w_s + S_\Gamma(w_y)](t,s,y)\} \geq 0.$$

With arguments similar to the ones used to show the subsolution property in [7], it can be proved that $w$ is also a subsolution of the equation

$$(4.13) \quad G(s, y, w_t(t,s,y), Dw(t,s,y), D^2w(t,s,y)) = 0.$$

We omit this proof because it has nothing to do with the small time behavior of stochastic integrals.

4.3. *DPE for the value function $v$.* For the derivation of the DPE for $v$ we have to restrict the control processes further by requiring that $\gamma$ is right-continuous and for all $t \in [0,T]$, there exists an $\varepsilon > 0$ such that almost surely,

$$(4.14) \quad \int_0^r |\gamma(u+t) - \gamma(t)|^2 \, du = O(r^{1+\varepsilon}) \quad \text{for } r \searrow 0.$$

Again, we proceed in three steps.

*Step* 1: *Dynamic programming.* We assume that for each $(t,s) \in [0,T) \times (0,\infty)$, there is an admissible control $(y,\nu) = (y,\alpha,\gamma)$ such that

$$X^\nu_{t,s,x,y}(T) \geq g(S_{t,s}(T)) \quad \text{where } x = v(t,s).$$

For a $(t,T]$-valued $\mathbb{F}^W$-stopping time $\tau$ and $\delta > 0$, we set $\tau_\delta := \tau \wedge (t+\delta)$. As in Section 4.2, it can be shown that

$$(4.15) \quad v(t,s) + \int_t^{\tau_\delta} Y^\nu_{t,s,y}(r)\,dS_{t,s}(r) \geq v(\tau_\delta, S_{t,s}(\tau_\delta)).$$

*Step* 2: *Application of Itô's lemma.* Again, we assume that the value function $v$ is smooth. Then, we can twice apply Itô's lemma in (4.15) to get for all $\delta > 0$:

$$-\int_t^{\tau_\delta} \xi(r)\,dr - \int_t^{\tau_\delta} \left\{ \zeta + \int_t^r \phi(u)\,du + \int_t^r \psi(u)\,dS_{t,s}(u) \right\} dS_{t,s}(r) \geq 0,$$

(4.16)

where

$$\xi(r) := \mathcal{L}v(r \wedge \tau, S_{t,s}(r \wedge \tau)),$$
$$\zeta := v_s(t,s) - y,$$
$$\phi(r) := \mathcal{L}v_s(r \wedge \tau, S_{t,s}(r \wedge \tau)) - \alpha(r),$$
$$\psi(r) := v_{ss}(r \wedge \tau, S_{t,s}(r \wedge \tau)) - \gamma(r),$$

and $\mathcal{L}$ is the Dynkin operator associated to the process $S$:

$$\mathcal{L}v(t,s) := v_t(t,s) + \tfrac{1}{2}\sigma_m^2 s^2 v_{ss}(t,s).$$

If we set

$$\tau := \inf\{r \geq t : |\log S_{t,s}(r) - \log s| > K\},$$

for some constant $K > 0$, then the processes $\xi$, $\phi$ and $\psi$ are bounded.

*Step* 3: *Small time path behavior of double stochastic integrals.* It follows from the boundedness of $\xi$ that there exists a constant $C_1 > 0$ such that for all $\delta > 0$,

(4.17)
$$\left| \int_t^{\tau_\delta} \xi(r)\,dr \right| \leq C_1 \delta.$$

From the boundedness of $\phi$ and Proposition 3.9 it can be deduced that

(4.18)
$$\lim_{\delta \searrow 0} \frac{1}{\delta} \left| \int_t^{\tau_\delta} \int_t^r \phi(u)\,du\,dS_{t,s}(r) \right|$$
$$= \lim_{\delta \searrow 0} \frac{1}{\delta} \left| \int_t^{t+\delta} \int_t^r \phi(u)\,du\, S_{t,s}(r)\sigma\,dW(r) \right| = 0.$$

Furthermore, since almost all paths of $S_{t,s}$ are Hölder-continuous of order $1/3$, it follows from Corollary 3.7(a) that

(4.19)
$$\limsup_{\delta \searrow 0} \frac{1}{h(\delta)} \left| \int_t^{\tau_\delta} \int_t^r \psi(u)\,dS_{t,s}(u)\,dS_{t,s}(r) \right|$$
$$= \limsup_{\delta \searrow 0} \frac{1}{h(\delta)} \left| \int_t^{t+\delta} \int_t^r \psi(u) S_{t,s}(u)\sigma\,dW(u)\, S_{t,s}(r)\sigma\,dW(r) \right|$$
$$< \infty.$$

It can be seen from (4.16) together with (4.17)–(4.19) that

$$\limsup_{\delta \searrow 0} \frac{1}{\sqrt{h(\delta)}} \int_t^{\tau_\delta} \zeta \, dS_{t,s}(r) \leq 0,$$

from which it can be derived by the Dambis–Dubins–Schwarz theorem and the law of the iterated logarithm for Brownian motion that $\zeta = 0$. Therefore, (4.16), (4.17) and (4.18) imply that

$$(4.20) \qquad \limsup_{\delta \searrow 0} \frac{1}{h(\delta)} \int_t^{\tau_\delta} \int_t^r \psi(u) \, dS_{t,s}(u) \, dS_{t,s}(r) \leq 0.$$

Since $\psi$ is right-continuous, it follows from (4.20) and Corollary 3.7(b) that $\psi(t) \leq 0$. Note that by the definition of $\psi$ and the gamma constraint (4.4),

$$(4.21) \qquad \Gamma_* \leq s^2(v_{ss}(t,s) - \psi(t)) \leq \Gamma^*.$$

By the boundedness and continuity of $\xi$, we obtain from (4.16) and (4.18) that

$$(4.22) \qquad \xi(t) \leq \liminf_{\delta \searrow 0} \frac{1}{\delta} \int_t^{\tau_\delta} \int_t^r -\psi(u) \, dS_{t,s}(u) \, dS_{t,s}(r).$$

Since $v$ is smooth, the process $v_{ss}(r, S_{t,s}(r))$ is almost surely locally Hölder-continuous of order $1/3$. Hence, since $\gamma$ satisfies (4.14), the process $\psi$ satisfies (4.14) as well. Therefore, we can apply Corollary 3.8(a) to conclude that

$$\liminf_{\delta \searrow 0} \frac{1}{\delta} \int_t^{\tau_\delta} \int_t^r -\psi(u) \, dS_{t,s}(u) \, dS_{t,s}(r)$$

$$= \liminf_{\delta \searrow 0} \frac{1}{\delta} \int_t^{t+\delta} \int_t^r -\psi(u) \, S_{t,s}(u)\sigma \, dW(u) \, S_{t,s}(r)\sigma \, dW(r)$$

$$= \frac{1}{2}\sigma^2 s^2 \psi(t),$$

which together with (4.22) shows that

$$\xi(t) \leq \tfrac{1}{2}\sigma^2 s^2 \psi(t),$$

that is,

$$(4.23) \qquad -v_t(t,s) - \tfrac{1}{2}\sigma^2 s^2 (v_{ss}(t,s) - \psi(t)) \geq 0.$$

Combined, (4.21) and (4.23) yield the following:

$$\hat{F}(v_t(t,s), s^2 v_{ss}(t,s)) := \sup_{\beta \geq 0} F(v_t(t,s), s^2 v_{ss}(t,s) + \beta) \geq 0,$$

where

$$F(p, A) := \min\{-p - \tfrac{1}{2}\sigma^2 A; \Gamma^* - A; A - \Gamma_*\}.$$

In [2] it is proved under weaker assumptions and with more general control processes that the value function $v$ is a viscosity solution of the equation

$$(4.24) \qquad \hat{F}(v_t(t,s), s^2 v_{ss}(t,s)) = 0.$$

4.4. *Discussion of the assumptions and related literature.* Using approximation arguments, it can be shown that $w$ is a viscosity solution of the DPE (4.13) without the assumption that it is always a minimum and smooth (see [7] for more details). Under additional continuity conditions on $\gamma$ it can also be shown that $v$ is a viscosity supersolution of the DPE (4.24) without assuming that it is always a minimum and smooth. For instance, with the arguments of Section 5 of [2] it can be shown without assumptions on $v$ that it is a viscosity supersolution of (4.24) if $\gamma$ is required to be of the form

$$(4.25) \qquad \gamma(r) = z + \int_t^r \gamma^1(u)\, du + \int_t^r \gamma^2(u)\, dW(u), \qquad r \in [t, T],$$

for $z \in \mathbb{R}$ and $\gamma^1, \gamma^2$ progressively measurable processes, and suitable boundedness conditions are satisfied. It is an open problem whether the supersolution property of $v$ can be shown without continuity assumptions on $\gamma$ such as (4.25) or (4.14). Another open problem is whether the value function $v$ corresponding to trading strategies of the form (4.3) with $\gamma$ of the form (4.25) is also a subsolution of (4.24).

However, assume that $g$ is continuous and let $v^{\mathrm{BS}}(t, s) = \mathrm{E}[g(S_{t,s}(T))]$ be the Black–Scholes price of $g(S_{t,s}(T))$. Then it follows from the comparison result, Proposition 3.9, in [2] that $v \geq v^{\mathrm{BS}}$ on $[0, T) \times (0, \infty)$, and $v > v^{\mathrm{BS}}$ on $[0, T) \times (0, \infty)$ whenever the function $g(s) + \Gamma^* \log s$ is not concave. On the other hand, without boundedness assumptions on the process $\alpha$ in (4.3), it follows from Theorem 4.4 in [1] that $v \leq v^{\mathrm{BS}}$ on $[0, T) \times (0, \infty)$ irrespective of the form of $g$.

To allow for a proof of a partial dynamic programming principle that is needed in the proof of the subsolution property of the value function, the control processes in [2] are also permitted to contain finitely many jumps. More precisely, the trading strategies in [2] are of the form

$$(4.26) \qquad Y(r) = \sum_{n=0}^{N-1} y^n \mathbf{1}_{\{\tau_n \leq r < \tau_{n+1}\}} + \int_t^r \alpha(u)\, du + \int_t^r \gamma(u)\, dS_{t,s}(u),$$
$$r \in [t, T],$$

where $t = \tau_0 \leq \tau_1 \leq \cdots$ is an increasing sequence of $[t, T]$-valued $\mathbb{F}^W$-stopping times such that the random variable $N := \inf\{n \in \mathbb{N} : \tau_n = T\}$ is bounded, all $y^n$ are $\mathcal{F}^W(\tau_n)$-measurable random variables and $\alpha, \gamma$ are $\mathbb{F}^W$-progressively measurable processes satisfying certain boundedness and continuity conditions (see Section 2.2 in [2]).

Denote by $v^{\mathrm{jumps}}$ the value function corresponding to this class of trading strategies. It is shown in [2] that $v^{\mathrm{jumps}}$ is the unique viscosity solution of (4.24) in a certain class of functions. Again, for continuous $g$, it follows from the comparison result, Proposition 3.9 in [2], that $v^{\mathrm{jumps}} \geq v^{\mathrm{BS}}$ on $[0, T) \times (0, \infty)$, and $v^{\mathrm{jumps}} > v^{\mathrm{BS}}$ on $[0, T) \times (0, \infty)$ whenever $g + \Gamma^* \log(s)$ is not concave. On the other hand, if the number of jumps $N$ in (4.26) is only required to be finite but not bounded, then it follows from Lemma A.3 in [3] that $v^{\mathrm{jumps}} \leq v^{\mathrm{BS}}$ on $[0, T) \times (0, \infty)$ for all $g$.

**Acknowledgments.** We are grateful to Peter Bank for pointing out the approximation results of [1] and [3] and to the anonymous referees for many useful suggestions. P. Cheridito gratefully acknowledges the kind hospitality of the Centre de Recherche en Economie et Statistique.

P. CHERIDITO
ORFE
PRINCETON UNIVERSITY
PRINCETON, NEW JERSEY
USA
E-MAIL: dito@princeton.edu

H. M. SONER
KOÇ UNIVERSITY
ISTANBUL
TURKEY
E-MAIL: msoner@ku.edu.tr

N. TOUZI
CENTRE DE RECHERCHE
EN ECONOMIE ET STATISTIQUE
PARIS
FRANCE
E-MAIL: touzi@ensae.fr
URL: www.crest.fr/pageperso/lfa/touzi/touzi.htm